\documentclass[a4paper,10pt]{article}

\usepackage[fleqn]{amsmath}
\usepackage{alltt}
\usepackage{amssymb}
\newtheorem{defi}{Definition}
\newtheorem{thm}{Theorem}

\author{Carl Bracken \\
School of Mathematical Sciences\\
University College Dublin\\
Ireland}

\newcommand{\done}{$\Box$\bigskip}

\begin{document}

%\begin{frontmatter}

%\date{\today}

\title{Pseudo Quasi-3 Designs and their Applications to Coding Theory}

\maketitle

%opening

\begin{abstract}
\noindent
We define a pseudo quasi-3 design as a symmetric design with the property that the derived and residual designs with respect to at least one block are quasi-symmetric. Quasi-symmetric designs can be used to construct optimal self complementary codes. In this article we give a construction of an infinite family of pseudo quasi-3 designs whose residual designs allow us to construct a family of codes with a new parameter set that meet the Grey Rankin bound.
\end{abstract}

\section{Introduction}

\noindent
Designs  were first considered for the purpose of designing statistical experiments, but 
have since found applications in other areas of mathematics. The study of quasi-symmetric designs began with S.S. Shrikhande \cite{Sh} who considered the duals of such designs. It was shown by  McGuire \cite{MG} that the existence of certain optimal error correcting codes was equivalent to the existence of particular quasi symmetric designs. One method of obtaining these quasi symmetric designs is by taking the derived and residual designs of a quasi-3 design (introduced by Cameron \cite{C}). Each quasi-3 design will give a pair of quasi symmetric designs and hence two error correcting codes. However, quasi-3 designs are quite rare so as an alternative approach we relax the conditions on the design so that it is not necessarily quasi-3 but can still give us the required quasi symmetric designs. We refer to such a design as pseudo quasi-3 and in this article we give a construction of an infinite family of these designs. This will allow us to obtain new quasi symmetric designs and their corresponding error correcting codes. We begin with some formal definitions.

\begin{defi} A t-design with parameters $t-(v,k,\lambda)$ is a pair $D(\mathcal{X},\mathcal{B})$ where $\mathcal{X}$ is a set of points of
cardinality $v$, and $\mathcal{B}$ a collection of k-element subsets of $\mathcal{X}$ called blocks, with the property that any $t$
points in $\mathcal{X} $ are contained in precisely $\lambda$ blocks.
\end{defi}
Throughout this article we have $t=2$, that is, we are only considering different types of $2$-designs.

\begin{defi} Let  $\mathcal{B} = \{B_1, B_2, \hdots B_b\}$ be the block set of a $t$-design and let $\mathcal{X} = \{X_1, X_2, \hdots X_v\}$ be its point set.
 . Then the incidence matrix $(M)$ of this design is the $b \times v$ binary array, with rows indexed by the blocks and columns
  indexed by the points of the design and the entry $M_{ij} = 1 $ if $ X_j \in B_i $ and $ 0 $  if $ X_j \notin B_i$, i.e.,
  the rows of $M$ are the characteristic vectors of the blocks as subsets of $\mathcal{X}$.
 \end{defi}

 \begin{defi} The dual design (denoted $D^T$) of a symmetric design $D$, is obtained by interchanging the point and block sets, while
   changing the relationship ``contained in" to ``contains".\end{defi}

\begin{defi} A symmetric design is a $2$-design, where the number of points equals the number of blocks.
\end{defi}

\noindent
Symmetric designs are sometimes called square designs as the incidence matrix is a square (not neccesarily a symmetric) matrix.
 It is shown in \cite{CvL} that the dual of a symmetric design is always a symmetric design.
 If $D$ is a symmetric design and $D^T$ is its dual, then the incidence matrix of $D^T$ will be the transpose
of the incidence matrix of $D$.

\begin{defi} A symmetric design is said to be quasi-3 (for blocks) if it has exactly two distinct triple block intersection sizes,
usually denoted $x$ and $y$ (with $x<y$).
 \end{defi}

\noindent
 A symmetric design can also be defined as quasi-3 for points, if there are only two possible numbers of blocks that contain any three
 points. It is clear that the dual of a design that is quasi-3 for points will be quasi-3 for blocks. The dual of a quasi-3 design is not necessarily quasi-3, see \cite{BR} for an example. Throughout this article the term
  ``quasi-3" shall denote quasi-3 for blocks.

\begin{defi}
A quasi symmetric design is a $2$-design with only two possible block intersection sizes.
\end{defi}

\begin{defi}
Let $D(\mathcal{X},\mathcal{B})$ be a $2-(v,k,\lambda)$ symmetric design, and $B$ a block.
 The derived design $(D_B)$ of $D(\mathcal{X},\mathcal{B})$
with respect to the block $B$ has block set $\mathcal{B} \setminus \{B\}$ and point set
 $\{x \in \mathcal{X} \ :\ x \in B \}$.
$(D_B)$ is a $2-(k,\lambda,\lambda -1)$ design.
\end{defi}
\begin{defi}
Let $D(\mathcal{X},\mathcal{B})$ be a $2-(v,k,\lambda )$ symmetric design, and $B$ a block.
 The residual design $(D^B)$ of $D(\mathcal{X},\mathcal{B})$
with respect to the block $B$ has block set $\mathcal{B} \setminus \{B\}$ and point set
 $\{x \in \mathcal{X} \ :\ x \notin B \}$.
$(D^B)$ is a $2-(v-k,k- \lambda ,\lambda )$ design.
\end{defi}

\begin{defi} A $u \times u$ Hadamard matrix $H$ is an $u \times u$ matrix with entries $1$ and $-1$, such that $$HH^{T}=uI.$$\end{defi}

\noindent
Hadamard matrices can only exist when $u$ is divisible by $4$ and are conjectured to exist for all such $u$. This has been verified for $u \leq 664$.

\begin{defi}
A normalised Hadamard matrix is a Hadamard matrix where both the first row and the first column consist entirely of ones.
\end{defi}

\noindent
As multiplying any row or column of a Hadamard matrix by minus one will retain the Hadamard matrix property,
as will permuting rows and columns, for any Hadamard matrix
there exists an equivalent normalised Hadamard matrix.

\section{Pseudo Quasi-3 Designs}
\noindent
The existence of a quasi-3 design with parameters $2-(4u^2, 2u^2-u, u^2-u)$ and triple block
intersection sizes of $ u^2/2-u $ and $ u^2/2-u/2$ implies the existence of two quasi-symmetric designs, taken as the derived and residual designs of the quasi-3 design. One with parameters $2-( 2u^2-u, u^2-u, u^2-u-1)$ and double block intersection sizes of $ u^2/2-u $ and $ u^2/2-u/2$, the other with parameters $2-(2u^2+u, u^2, u^2-u)$ and block intersection sizes $u^2/2-u/2$ and $u^2/2$. When $u$ is a power of $2$ these quasi-3 designs can be constructed (see \cite{B} and \cite{JT2}) and hence the quasi symmetric designs are obtained. However, when $u$ is not a power of two the existence of these designs is an open problem. This family of quasi-3 designs could exist for all even $u$, no case has been ruled out or proven. It is possible that these quasi symmetric designs exist even if the quasi-3 design does not. The purpose of this article is to consider symmetric designs with a weaker property than the quasi-3 property which still give a pair of quasi symmetric designs as derived and residual designs. We are thus motivated to define the following.

\begin{defi}
A pseudo quasi-3 design is a symmetric $2$-design with the property that the block intersection sizes of all triples of blocks that contain one specified block, takes one of two distinct values.
\end{defi}

\noindent
As a convention we shall place the specified block as the first row of the incidence matrix of the design. This means that any triple of rows that contains the first row must have one of two possible intersection numbers. It should be noted that the design may have many blocks that could have this triple intersection property but in order for the design to be pseudo quasi-3, it only needs one. If all blocks are such that the triples containing them have only two possible intersection numbers then the design is quasi-3.
We now obtain two quasi symmetric designs as the derived and residual designs of a pseudo quasi-3 design by projection onto the first block. In the sequel we offer a construction of pseudo quasi-3 designs. In the final section we discuss the new families of optimal error correcting codes that can be constructed from the quasi symmetric designs obtained from the pseudo quasi-3 designs.

\section{Pseudo Quasi-3 Designs from Hadamard Matrices}
In this section we give a construction of pseudo quasi-3 designs with parameters $2-(4u^2, 2u^2-u, u^2-u)$ and triple block intersections of $ u^2/2-u $ and $ u^2/2-u/2$ provided the triple contains the first block, whenever there exist a $u \times u$ Hadamard matrix.

\noindent
{\bf The Construction.}
Let $H_u$ be any normalised $u \times u$ Hadamard matrix and let
$$H_{2u}=  \left(\begin{array}{cc} H_u & \ \ H_u \\ H_u & -H_u
 \end{array} \right).$$
$H_{2u}$ is also a normalised Hadamard matrix. Let $h_i$ denote the $i^{th}$ row of $H_{2u}$ and define
$\tilde{A_i}:=h_i \otimes (h_i)^T$ and $\tilde{S_i}:=h_u \otimes (h_i)^T$, for \ $0 \leq i \leq 2u-1$, where $\otimes$ denotes the Kronecker product. 
Let $L$ be the $2u-1 \times 2u-1$ Latin square defined on the symbols $0$ to $2u-1$ with $u$ omitted by taking the cyclic shifts of
 \begin{tabular}{|c|c|c|c|c|c|c|c|c|c|c|c|}
 \hline
 $0$ & $1$ & $2$ & $\hdots$ & $u-2$ & $u-1$ & $2u-1$ & $2u-2$ & $\hdots$ & $u+3$ & $u+2$ & $u+1$\\
 \hline
 \end{tabular}\\

\noindent
Next we place matrix $A_i$ in position $i$ in $L$ to obtain the matrix $L(A)$. Now we use $L(A)$ and the $S_i$ matrices, as well as their transposes (denoted $S_i^T$) to construct the following matrix which we name $P_u$.

\begin{center}
$P_u:=$
\begin{tabular}{|c|c|c|c|c|c|}
 \hline 
 $0$  & $S_1 $ & $ S_2 $ & $ S_3 $ & $ \hdots $ & $ S_{2u-1}$\\
 \hline
 $S_1^T$ & \multicolumn{5}{c|}{} \\ \cline{1-1}
 $S_2^T$ & \multicolumn{5}{c|}{} \\  \cline{1-1}
 $S_3^T$ & \multicolumn{5}{c|}{$L(A)$} \\ \cline{1-1}
$ \vdots$ &  \multicolumn{5}{c|}{} \\ \cline{1-1}

 $S_{2u-1}^T$ & \multicolumn{5}{c|}{} \\
\hline
 \end{tabular}\\
\end{center}
In the following theorem we demonstrate that the above matrix is a pseudo quasi-3 design.

\begin{thm}
Let $P_u$ be constructed as above. Then $P_u$ is the incidence matrix of a pseudo quasi-3 design with parameters
$2-(4u^2, 2u^2-u, u^2-u)$ and triple block intersections of $ u^2/2-u $ and $ u^2/2-u/2$ when the triple contains at least one block from the first $2u$ blocks.
\end{thm}
Proof:\\
\noindent
First we note that $A_0$ is the all zero matrix and the other $A_i$ matrices and the $S_i$ matrices have $2u$ rows and columns with each row having $u$ $1$'s and $u$ $0$'s. This establishes $v=4u^2$ and $k=2u^2-u$.

Any two rows of a particular $A_i$ or $S_i$ matrix, with $i \neq 0$, intersect in $u$ or $0$ places. If we compare two rows of $P_u$ from the same row of cells, the number of times they agree is determined by the agreement of the corresponding columns of $H_{2u}$.

As $H_{2u}$ is a Hadamard matrix, any two columns agree in $u$ places. However, one of these agreements corresponds to the agreement of any two rows of the all zero matrix. The remaining $u-1$ agreements yield intersections of $u$ in each cell. This gives any two such rows in $P_u$ an intersection of $u(u-1)=u^2-u$. If we compre two rows of $P_u$ from different rows of cells, we have two positions with no intersection due to the presence of exactly one all zero matrix in every row of cells. In the remaining $2u-2$ pairs of cells we have an intersection of $\frac{u}{2}$ as any two rows of $H_{2u}$ have $\frac{u}{2}$ positions in which both rows have $-1$'s. Therefore the intersection of any two such rows is $\frac{u}{2}(2u-2)=u^2-u$. This verifies that $P_u$ is a symmetric design with $\lambda = u^2-u$.

Next we consider the intersection of three blocks in which at least one of the blocks is from the first row of cells. If all three rows are from the first row of cells then the triple intersection is $u$ times the triple intersection of points in the Hadamard 3-design formed by $\left(\begin{array}{c} H_{2u}^* \\ H_{2u}^{*c}
\end{array}\right)$, where $H_{2u}^*$ is $H_{2u}$ with $1$ and $-1$ changed to $0$ and $1$ respectively and $H_{2u}^{*c}$ is the complement of $H_{2u}^*$. These designs have triple intersections of $u-1$, so the three rows of $P_u$ intersect in $\frac{u^2}{2}-u$ places.

If we take two rows from the first row of cells and one from another row, the intersection of the first two rows is in $u-1$ sections of length $u$ at the first or last $u$ positions of each cell. As the first or last $u$ positions of each $A_i$ cell, when $i$ is not $u$ or $0$, has $\frac{u}{2}$ $1$'s and $\frac{u}{2}$ $0$'s we get an intersection of $\frac{u}{2}(u-1)=\frac{u^2}{2}-\frac{u}{2}$ if the $A_0$ cell is not involved and 
 $\frac{u}{2}(u-1)-\frac{u}{2}=\frac{u^2}{2}-u$ if $A_0$ is involved in the intersection.

A similar argument applies when we consider one row from the first row of cells and two rows from some other row of cells.

Finally we consider the case when we take one row from the first row of cells and the other two rows from two different rows of cells. Recall the intersection of any two rows from distinct matrices $A_i$ and $A_j$, with $i$ and $j$ both non zero, is $\frac{u}{2}$. This comprises of an intersection of $\frac{u}{4}$ in the first $u$ positions and an intersection of $\frac{u}{4}$ in the last $u$ positions, unless $|i-j|=u$ in which case we have an intersection of $\frac{u}{2}$ in the first $u$ positions and an intersection of $0$ in the last $u$ positions or an intersection of $0$ in the first $u$ positions and an intersection of $\frac{u}{2}$ in the last $u$ positions.

We claim that, in the position-wise differences of any two rows of $L$, precisely one of the resulting differences has magnitude $u$. To demonstate this claim we reduce all entries in $L$ modulo $u$ and observe that the resulting array is the table of Lee differences for the elements of $\mathbb{Z}_{2u-1}$. As any two elements of  $\mathbb{Z}_{2u-1}$ have the same Lee distance to exactly one other element, we have zero appearing exactly once in the position-wise differences of any two rows of this array. Therefore a multiple of $u$ will appear precisely once in the differences of any two rows of the Latin square. With $\pm u$ being the only possibilities the claim is proven.

As the rows of any $S_i$ matrix consists of $u$ $0$'s followed by $u$ $1$'s or vise versa, to consider a triple intersection involving one of the first $2u$ blocks, we need only consider the intersection of the other two blocks when restricted to the first or last $u$ positions in each pair of cells. This gives a triple intersection consisting of three triples of cells with no intersection, one triple of cells with intersection of $0$ or $\frac{u}{2}$ and the remaining $2u-4$ triples intersecting in $\frac{u}{4}$ places. This yields triple intersections of $ u^2/2-u $ and $ u^2/2-u/2$ as required. 
 \done

We can now take the derived and residual designs with respect to any one of the first $2u$ blocks of $P_u$ to obtain two quasi symmetric designs. 

The derived design has parameters $2-(2u^2-u, u^2-u, u^2-u-1)$ and double block intersection sizes of $ u^2/2-u $ and $ u^2/2-u/2$.

The residual design has parameters $2-(2u^2+u, u^2, u^2-u)$ and double block intersection sizes $u^2/2-u/2$ and $u^2/2$.

From this construction we now know that these designs exist whenever there exists a $u \times u$ Hadamard Matrix, which is virtually every $u$ that is a multiple of four.
The parameter set for the first of these designs was already achieved in \cite{Bra} with a different construction. The parameters for the residual design are new.

\section{Application to Coding Theory}
An error correcting code $C$ with parameters $(n,M,d)$ over an alphabet $A$ is a subset of $A^n$ with the properties that $|C|=M$ and any two elements of $C$ differ in at least $d$ coordinates. The elements of $C$ are called words and when $A=\{0,1\}$, we say that $C$ is self-complementary if for all $c \in C$, $\bar{c} \in C$ where $\bar{c}$ denotes the binary complement of $c$.
 In \cite{MG} the existence of certain quasi-symmetric designs was shown to be equivalent to the existence of a self complementary code meeting the Grey-Rankin bound with equality. The Grey-Rankin bound is an upper bound on $M$ for a fixed $n$ and $d$. It states that
 $$  M \le \frac{8d(n-d)}{n-{(n-2d)}^2} $$
 for any $(n, M, d)$ self-complementary code provided the right hand side of the inequality is positive.
The following theorem is ``Theorem A (part (ii))"  from \cite{MG}.
\begin{thm}
 For $n$ even and $n- \sqrt{n}<2d<n $ there exists a self-complementary $(n,M,d)$ code with
 $M=8d(n-d)/(n-{(n-2d)}^2)$ if and only if $d$ is even and there exists a
quasi-symmetric $2-(n,d,\lambda ) $ design with block intersection sizes $d/2$ and $(3d-n)/2$, where
 $\lambda =d(d-1)/(n-{(n-2d)}^2).$
\end{thm}

In \cite{Bra} the infinite family of error correcting codes with parameters $(2u^2-u, 8u^2, u^2-u)$ were constructed using $u \times u$ Hadamard matrices. These parameters satisfy the Grey-Rankin bound with equality, therefore the above theorem implies the existence of quasi symmetric designs with parameters $2-( 2u^2-u, u^2-u, u^2-u-1)$. The $4u^2-1$ words of weight $u^2-u$ in this code form the incidence matrix of the quasi symmetric design.

The residual design of $P_u$ constructed in the last section has parameters $2-(2u^2+u, u^2, u^2-u)$. If we take its incidence matrix and its complement along with the all-zero and all-one words, we have a self complementary code with parameters $(2u^2+u, 8u^2, u^2)$. Again, these parameters meet the Grey-Rankin bound and hence the codes in this new family are optimal.

\section{Closing Remarks and Open Proplems}
In this article we have shown that there exists quasi symmetric designs with parameters
 $2-( 2u^2-u, u^2-u, u^2-u-1)$ and  $2-(2u^2+u, u^2, u^2-u)$ whenever there exists a $u \times u$ Hadamard matrix. However, these parameters are permissable for any even $u$, so when $u$ is not divisible by $4$ another construction is needed. The only examples of quasi symmetric designs with the above parameters and $u$ not a multiple of $4$ were constructed in \cite{BMW} with $u=6$. This pair of designs were not taken as derived and residual designs of a pseudo quasi-3 design and the existence of such a design is open.

{\bf Open Problem 1.} Does there exist a pseudo quasi-3 design with parameters $(144, 66, 30)$ and triple block intersections of $12$ and $15$?

It can be seen from the construction in \cite{BMW} that the non-existence of such a design implies the non-existence of the projective plane of order $12$.

Even if we cannot construct any more pseudo quasi-3 designs, it may be possible to obtain the remaining quasi symmetric designs by other constructions.

{\bf Open Problem 2.} Construct quasi symmetric designs with parameters
$2-( 2u^2-u, u^2-u, u^2-u-1)$ and  $2-(2u^2+u, u^2, u^2-u)$ when $u>6$ and not divisible by $4$.

\end{document}